\documentclass[12pt]{amsart}
\usepackage[english]{babel}
\usepackage{amsmath}
\usepackage{amssymb}
\usepackage{latexsym}
\usepackage{amscd}
\newtheorem{theorem}{Theorem}[section]         
\newtheorem{lemma}[theorem]{Lemma}             
\newtheorem{proposition}[theorem]{Proposition} 
\newtheorem{remark}[theorem]{Remark}           
\newtheorem{definition}[theorem]{Definition}   
\newtheorem{example}[theorem]{Example} 
\numberwithin{equation}{section}

\newenvironment{namelist}[1]{%
\begin{list}{}
    {
      
      \settowidth{\labelwidth}{#1}
      \setlength{\leftmargin}{1.1\labelwidth}
    }
  }{%
\end{list}}
\newcommand{\Fa}{\ensuremath{\mathcal{F}}}

\newcommand{\C}{\ensuremath{\mathcal {C}}}
\newcommand{\Pa}{\ensuremath{\mathcal {P}}}
\newcommand{\F}{\ensuremath{\mathcal {M}}}

\newcommand{\rk}{r}
\newcommand{\M}{\ensuremath{ M}}
\newcommand{\cl}{\rm{cl}}
\begin{document}
\title[chordal graphs and matroids]
{how is a chordal graph like a supersolvable binary matroid?}
\author{Raul Cordovil, David Forge and Sulamita Klein}
\address{{}\newline
Departamento de Matem\'atica,\newline 
Instituto Superior T\' ecnico \newline
 Av.~Rovisco Pais
 - 1049-001 Lisboa  - Portugal}
\email{cordovil@math.ist.utl.pt}
\thanks{ The  first  author's research was 
supported in part by FCT/FEDER/POCTI (Portugal)  and
 the project SAPIENS/FEDER/36563/00.  
 The third author  was partially supported by CNPq,
MCT/FINEP PRONEX Project 107/97, CAPES (Brazil)/COFECUB (France),
project number 213/97, FAPERJ}
\address{\newline 
LRI, UMR 8623, Batiment 490 Universit\'e Paris-Sud\newline
91405 Orsay Cedex, France}
\email{forge@lri.fr}
\address
{{}\newline
Instituto de Matem\'atica and COPPE-Sistemas,\newline
Universidade Federal do Rio de Janeiro,\newline
Caixa Postal 68511, 21945-970, Rio de janeiro, RJ, Brasil
}
\email{sula@cos.ufrj.br}
\dedicatory{To the memory of Claude Berge}
\begin{abstract}
{ Let $G$ be a  finite simple graph. From the pioneering work 
of R. P. Stanley it is known that the cycle matroid of
$G$ is supersolvable iff $G$ is chordal (rigid): 
this is another way to read Dirac's
theorem on chordal graphs. 
 Chordal binary matroids  are not in general
supersolvable. Nevertheless we prove that,  
for every supersolvable
binary matroid $M$, a
maximal chain of modular flats of $M$  canonically   determines  a
chordal graph. }
\end{abstract}
\maketitle 
\section{Introduction and notations} 
  Throughout this note  $\M$ denotes
 a  matroid of rank $r$ on the ground set
 $[n]:=\{1,2,\dotsc, n\}$.   We refer to 
\cite{Ox,W1}  
 as standard sources  for  matroid theory. We recall and fix
 some notation of matroid theory. The restriction of $\M$ 
to a subset $X\subseteq [n]$ is denoted $\M|X.$ A matroid $\M\,$ is said to be
\emph{simple} if all  circuits have at least three elements.  A matroid
$\M\,$ is  \emph{binary}  
if the symmetric difference of any two different 
circuits of $M$\, is a  union of disjoint  circuits. Graphic and 
cographic matroids are extremely important
examples of binary matroids.   
The dual of $\M$ is denoted $\M^*$. Let
$\C=\C(\M)$ [resp.
$\C^*=\C^*(\M)=\C(\M^*)$] be the set of  circuits 
[resp. cocircuits] of $\M.$ Let $\C_\ell:=\{C\in \C:~ |C|\le \ell\}.$ 
In the following the singleton 
$\{x\}$ is denoted by  $x$. We will denote by 
$$\cl(X):=X\cup\{x\in[n]:\,\exists C\in \C\,, 
C\setminus X=x\},$$    the  {\em closure}  in \M\,  of a subset
$X\subseteq [n].$  We say that $X\subseteq {[n]}$ is a
\emph{flat} of \M \, if $X=\cl(X).$
The set $\Fa(\M)$ of flats of \M, ordered by inclusion,
is a geometric lattice. The {\em rank} of a flat $F\in \Fa,$ 
denoted $r(F)$, is equal to $m$ if there are $m+1$ flats in
a maximal chain of flats  from
$\emptyset$ to $F$. The flats of rank 1, 2, 3 and 
$r-1$ are called {\em points, lines, planes,} and 
{\em hyperplanes}   respectively. A line $L$
with two  elements is called {\em  trivial} and a line with three
elements is called nontrivial (a binary matroid
has no line with more than three points).  Given a set
$X\subseteq [n]$, let  $\rk(X):=\rk(\cl(X))$. A pair
$F, F'$ of flats is called \emph{modular} if
$$ \rk(F)+ \rk(F')= \rk(F \vee F')+ \rk(F \wedge F'). $$
A flat $F\in\Fa$ is \emph{modular} if it forms a modular pair with every 
other flat $F'\in\Fa.$ The notion of supersolvable lattices
was introduced and studied by Stanley in  \cite{Stan2}. In the 
particular case of  geometric lattices the definition can be read as follows.
\begin{definition}\cite{Stan2}.
{\em A matroid $M$ on $[n]$ of rank $r$ is  {\em supersolvable} 
 if there is a 
maximal chain of modular flats  \F
$$\F:=\, F_0(=\emptyset)\subsetneq \dotsm\subsetneq F_{r-1}\subsetneq
F_r(=[n]).$$ We call  $\F$ an {\em $\M$-chain}\, of $\M$.
To the $\M$-chain $\F$ we associate the partition  $\Pa$ of $[n]$ 
 $$\Pa:=\, F_1\uplus\dotsm \uplus (F_i\setminus
F_{i-1})\uplus\dotsm
\uplus( F_r\setminus F_{r-1}).$$
We call $\Pa$ an {\em $\M$-partition} of $\M.$ 
}\end{definition}
 We recall that a graph $G$ is said {\em chordal} (or {\em rigid} or {\em
triangulated})  if every cycle of
length at least four has a chord. Chordal graphs are treated extensively
in Chapter 4 of \cite{col}. 
The
notion of a ``chordal matroid"  has also been recently explored  in
the literature, see
\cite{B}.  
\begin{definition}[\cite{Ba} p. 53]\label{rigid}
{\em Let $M$ be an arbitrary matroid (not necessarily simple or binary).
A circuit $C$  of 
$M$ has a {\em chord}  
$\,i_0$ if there are two circuits $C_1$ and $C_2$
such that $C_1\cap C_2=i_0$ and $C=C_1\Delta C_2.$ In this case, we say that 
the chord $i_0$ {\em splits} the circuit $C$ into the circuits $C_1$ and $C_2.$
We say that a  matroid is $\ell$-\emph{chordal}  if  
every circuit with at least $\ell$ elements has a chord. 
A simple matroid $\M$ is  \emph{chordal} if
it is 4-chordal. 
}
\end{definition}
In this paper we  always suppose that the edges of a graph $G$ 
are labelled with the integers of $[n]$.
If nothing else is indicated we suppose $G$ is a connected graph.
Let $\M(G)$ be the \emph{cycle matroid} of the graph 
$G$: i.e., the elementary cycles of $G,$ as subsets of $[n],$ are  
the circuits of $\M(G).$ In the same way, the minimal cutsets 
of a connected  graph $G$ (i.e, a set of edges that 
disconnect the graph) are the
circuits of a matroid on $[n],$  called the
\emph{cocycle matroid} of $G.$ A matroid is  \emph{graphic} (resp. 
\emph{cographic}) if 
it is the cycle (resp. cocycle) matroid of a graph. 
The cocycle matroid of $G$ is dual to the 
cycle matroid of $G$ and  both are binary.
The cocycle matroids of the complete graph $K_5$ and of the complete
bipartite graph $K_{3,3}$ are examples of binary but not graphic matroids; see
Section 13.3 in \cite{Ox} for details. The Fano matroid  
is an example of a supersolvable binary 
matroid that is neither graphic nor cographic. Finally, note 
that an elementary cycle
$C$ of $G$ has a chord iff $C$ seen as a circuit of the matroid
$\M(G)$ has a chord.
\begin{example}\label{exa2}
{\em
 Consider the chordal graph $G_0=G_0(V,[7])$ in  Figure~1 and 
the corresponding cycle matroid $\M(G_0)$. It is clear that
$$\F:= \emptyset  \subsetneq\{1\} \subsetneq \{1,2,3\} 
\subsetneq \{1,2,3,4,5\} \subsetneq [7]$$
is an $M$-chain. The associated $M$-partition  is
$$\Pa:= \{1\}\uplus \{2,3\}\uplus \{4,5\}\uplus \{6,7\}.$$
The linear order of the vertices is such that for every 
$i$ in $ \{2,3, 4, 5\}$ the neighbors of the vertex $v_i$ contained in the set
$\{v_1,\dotsc v_{i-1}\}$ form a clique; this is  Dirac's well known 
characterization of chordal graphs  (see
\cite{Dirac, col}). This is also a characterization of graphic
supersolvable matroids (see Proposition 2.8 in \cite{Stan2}).
That is, a graphic matroid $\M(G)$ is supersolvable iff\, 
the vertices of $G$ can be labeled as $v_1,v_2,\dotsc,
v_m$ such that, for every $i=2,\dotsc, m,$ the  neighbors of 
$v_i$ contained in the
set
$\{v_1,\dotsc v_{i-1}\}$ form a clique.
We say that a linear order of the vertices of $G$ with the above properties 
is an {\em S-label} of the vertices  of $G.$ }
\end{example}
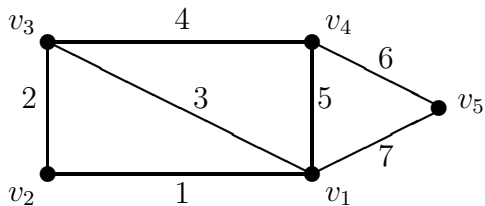
\begin{figure}
\thicklines
\begin{center}
\begin{picture}(0,130)(80,0)
\put(10,60) {\line(2,0){100}}
\put(10,110) {\line(2,0){100}}
\put(10,110) {\line(2,-1){100}}
\put(10,60) {\line(0,1){50}}
\put(10,60) {\line(2,0){100}}
\put(110,60) {\line(0,1){50}}
\put(-5,115){$v_3$}
\put(115,115){$v_4$}
\put(115,50){$v_1$}
\put(-5,50){$v_2$}
\put(58,115){4}
\put(58,49){1}
\put(0,85){2}
\put(65,85){3}
\put(112,85){5}
\put(135,100){6}
\put(135,63){7}
\put(158,85){\circle*{6}}
\put(165,85){$v_5$}
\put(10,60) {\circle*{6}}
\put(10,110) {\circle*{6}}
\put(110,60) {\circle*{6}}
\put(110,110) {\circle*{6}}
\put(110,110){\line(2,-1){48}}
\put(110,60){\line(2,1){48}}
\end{picture}
\vspace{-15mm}
\caption{Graph $G_0$}
\end{center}
\end{figure}
Ziegler proved  that every supersolvable 
binary matroid without a Fano submatroid is graphic  (Theorem 2.7 in \cite{Z}).
In the next section we answer the following natural question:
\begin{namelist}{xxx}
\item[~$\circ$]{\em For a generic binary
matroid, what are the relations between the notions of 
``chordal" and ``supersolvable" ?}
\end{namelist}
\section{chordal and supersolvable matroids}
\begin{lemma}\label{chord2} Let $\M$ be a simple binary matroid. 
The following two conditions are equivalent for every
circuit
$C$ of
$\M$:
\begin{namelist}{xxxxxx}
\item[$(\ref{chord2}.1)$] $C \subsetneq\cl(C),$
\item[$(\ref{chord2}.2)$] $C$ has a chord.
\end{namelist}
For  nonbinary matroids only the implication
$(\ref{chord2}.2)\Rightarrow (\ref{chord2}.1)$  holds.
\end{lemma}
\begin{proof}
If $i\in \cl(C)\setminus C$, then there is a circuit $D$ such 
that $i \in D$ and $D\setminus i\subsetneq C$. As $\M$ is binary
$D'=D\Delta C$ is also a circuit of $\M.$ So $i$ is a chord of 
$C$. If $i$ is a chord of $C$, then clearly  $i\in \cl(C)$.
Finally, in the uniform rank-two nonbinary matroid $U_{2,4}$, the set
$C=\{1,2,3\}$ is a   circuit without a chord  but $C \subsetneq\cl(C)=[4]$.
\end{proof}
\begin{theorem} A binary supersolvable matroid  $M$ is chordal but  
the converse does not hold in general. 
\end{theorem}
\begin{proof}
Let $\F:=\,\,\emptyset\subsetneq \dotsm\subsetneq F_{r-1}\subsetneq F_r=[n]$
be an $M$-chain of \M.
Suppose by induction that the restriction of $\M$ to $F_{r-1}$ is chordal. 
The result is clear in the case that 
 $C^*:=[n]\setminus F_{r-1}$ 
 is a singleton.
  Suppose that $|C^*|>1$   and consider a circuit 
$C$ of $\M$ not
contained in the modular hyperplane $F_{r-1}.$ Then   there are two 
elements $i,j\in C\cap C^*$ and the line
$\cl(\{i,j\})$ meets 
$F_{r-1}.$ So $C \subsetneq\cl(C)$ and we know  from Lemma~\ref{chord2}  
that C has a chord.\newline
A counterexample of the converse is $\M^*(K_{3,3})$, the cocycle  matroid of
the complete bipartite graph
$K_{3,3}$. It is easy to see from its geometric representation that it
is chordal but not supersolvable (see \cite{Z} and page 514 in \cite{Ox} for
its geometric representation). 
\end{proof}
\begin{definition}[\cite{CR}]\label{ell} 
{\em Let $\M$ be an arbitrary matroid and consider an integer
$\ell\geq 2.$  The matroid $\M$ is {\em $\ell$-closed}  if
the following two conditions are equivalent for every subset $X\subseteq [n]:$
\begin{namelist}{xxxxxx}
\item[$(\ref{ell}.1)$] $X$ is closed,
 \item[$(\ref{ell}.2)$] for every subset $Y$ of $X$ with at 
most $\ell$ elements we have $\cl(Y)\subseteq X.$
\end{namelist}
}
\end{definition}
We note that Condition $(\ref{ell}.2)$ is equivalent to:
\begin{namelist}{xxxxxx}
\item[$(\ref{ell}.2')$] for every circuit $C$ of $\M$  with at most 
$\ell+1$ elements $$|C\cap X|\geq |C|-1
\Longrightarrow C\subseteq X.$$
\end{namelist}
\begin{definition}\label{Db}
{\em Let $\C'$ be a subset of $\C$, the set of circuits of $\M.$ 
Let $\cl_{\Delta}(\C')$ denote the smallest subset of $\C$ such that:
\begin{namelist}{xxxxxxx}
\item[$(\ref{Db}.1)$] $\C'\subseteq \cl_{\Delta}(\C'),$
\item[$(\ref{Db}.2)$] whenever a circuit splits into two circuits $C_1$ and
$C_2$ that are in $\cl_{\Delta}(\C')$ then $C$ is also in $\cl_{\Delta}(\C').$
\end{namelist}
}
\end{definition}
\begin{theorem}\label{main} For every simple  binary matroid $\M$ the 
following three conditions are equivalent:
\begin{namelist}{xxxxxxxx}
\item[$(\ref{main}.1)$] $\M$ is $\ell$-closed,
\item[$(\ref{main}.2)$] $\M$ is $(\ell+2)$-chordal,
\item[$(\ref{main}.3)$] $\C(\M)=\cl_\Delta(\C_{\ell+1}).$
\end{namelist}
\end{theorem}
\begin{proof}
$(\ref{main}.2)\Longleftrightarrow (\ref{main}.3).$ This equivalence is a 
direct consequence of the definitions.\newline
$(\ref{main}.1)\Longrightarrow (\ref{main}.2).$ Consider a circuit $C$
with at least
$\ell+2$ elements and suppose for a contradiction that $C$ is not  chordal.
 From Lemma~\ref{chord2} we know that $\cl(C)=C.$ Pick an element $i\in C.$
Then the set $X=C\setminus i$ is not closed but every subset $Y$ of $X$ with 
at most $\ell$ elements is closed which is a contradiction.\newline 
$(\ref{main}.3)\Longrightarrow (\ref{main}.1)$.  Let $X$ be a subset of $[n]$ 
and suppose that for every  circuit $C$   with at most $\ell+1$ elements
such that $|C\cap X|\geq |C|-1,$ we have $C\subseteq X$; see $(\ref{ell}.2')$.
 To prove that $X$ is closed it is enough to prove that for every 
circuit $C$ such that $|C\cap X|\geq |C|-1,$ we have $C\subseteq X.$ 
Suppose that the result is
true for every circuit with at most $m$ elements and let $D$ be a 
circuit with $m+1$ elements such that $D\setminus d\subset  X$ with 
$d\in D.$ By hypothesis
there are circuits
$C_1, C_2\in \cl_\Delta(\C_{\ell+1})$ such that $C_1\cap C_2=i$
and $D=C_1\Delta C_2.$  Suppose  w.l.o.g that $d\in C_1$.
We have $C_2\setminus i\subset X$ and since $|C_2|<|D|$ we also have
$i\in C_2\subset X$. We have now that $C_1\setminus d\subset X$ and 
$|C_1|\le m$ we also have
that $C_1\subset X$. This gives that $D\subseteq X$ and concludes the proof.
\end{proof}
We make use of the following elementary but useful proposition which 
is a particular case of Proposition 3.2 in
\cite{Stan2}. The reader can easily check it from 
Brylawski's characterisation of modular hyperplanes.
 \begin{proposition}
Let   $M$ be a supersolvable matroid and  $$\F:=\,F_0\subsetneq 
\dotsm\subsetneq F_{r-1}\subsetneq F_r$$ an $M$-chain.  Let $F$ be a flat
of $M$. Then $M|F$, the restriction of $M$ to the flat $F$, 
is a supersolvable matroid and $\{F_i\cap F: F_i\in 
\F\}$ is the set of (modular) flats of  an\,\, $M|F$-chain.\qed 
\end{proposition}
\begin{definition}\label{Pa}
{\em Let   $\Pa= P_1\uplus\dotsm \uplus P_r$ be an $M$-partition of 
a supersolvable matroid $\M$.
We associate to  $(\M, \Pa)$ a graph $G_{\Pa}$ such that:
\begin{namelist}{xxxxxx}
\item[~$\circ$]  $V(G_{\Pa})=\{P_i: i=1,2,\dotsc,r\}$ is the vertex 
set of $G_{\Pa},$
\item[~$\circ$] $\{P_i,P_j\}$ is an edge of $G_\Pa$ iff 
there is a nontrivial line $L$ of $\M$ meeting $P_i$ and $P_j$.
 \end{namelist}
We call  $G_{\Pa}$ the {\em S-graph} of the pair $(\M, \Pa)$.}
\end{definition}
Note that every nontrivial line $L$ of the binary supersolvable 
matroid $M$ meets exactly two
$P_i's$ and if $L$ meets $P_i$ and $ P_j$, with $i<j$, 
necessarily $|P_i\cap L|=1$ and
$|P_j\cap L|=2$. Indeed $F_{j-1}=\bigcup_{\ell=1}^{j-1} P_\ell$ is 
a modular flat disjoint from $P_j$, so $|F_{j-1}\cap L|=1$. 
This simple property
will be used  extensively  in the proof of Theorem \ref{chord}.
Given a chordal graph $G$ with a fixed S-labeling, we get an
associated supersolvable matroid $M(G)$ and an associated $M$-partition $\Pa.$
We say that $G_\Pa,$ the S-graph determined by  $(\M(G),\Pa),$ 
is the {\em derived S-graph} of $G$ for this S-labeling.
\begin{remark}\label{sub}
{\em 
Note that the derived S-graph $G_\Pa$ of a chordal graph $G$ is a 
subgraph of $G.$
Indeed set $V(G_\Pa)=\{P_1,\dotsc, P_m\}$ and consider the map 
$P_\ell\to v_{\ell+1},\, \ell=1,\dotsc m.$ 
Let $\{P_i, P_j\},\, 1\leq i<j\leq m$, be an edge of $G_\Pa.$ 
From the definitions we see that
$\{v_{i+1}, v_{j+1}\}$ is necessarily an edge of $G.$
}
\end{remark} 
\begin{example}
{\em Consider the S-labeling of the graph $G_0$ given in Figure~1 and 
the associated  $M$-partition $\Pa$ (see 
Example~\ref{exa2}). The derived  S-graph $G_{\Pa}$ is a path from $P_1$ 
to $P_4$.
Consider now the $M$-partition of $\M(G_0):$
$$\Pa':=  \{4\}\uplus \{3,5\}\uplus \{1,2\} \uplus \{6,7\}.$$
In this case the corresponding S-graph $G'_{\Pa'}$ is $K_{1,3}$ with
$P_2$ being the degree-3 vertex.
It is easy to prove that for any $M$-partition $\Pa$ of the cycle matroid
of the complete graph
$K_\ell,$  the S-graph $G_{\Pa}$ is the complete graph $K_{\ell-1}.$}
\end{example}
Our main result is:
\begin{theorem}\label{chord} Let $\M$ be a simple binary 
supersolvable matroid with an $M$-partition $\Pa$. 
Then the S-graph $G_\Pa$ is  chordal. 
\end{theorem}
\begin{figure}
\begin{center}
\begin{picture}(-140,230)(100,10)
\put(60,220) {\circle*{6}}
\put(52,223) {$x$}
\put(72,188) {\circle*{6}}
\put(77,188){$x'$}
\put(60,160){\circle*{6}}
\put(65,160){$y'$}
\put(14,130){\circle*{6}}
\put(3,130){$y$}
\put(60,160) {\line(-3,-2){130}}
\put(102,98){\circle*{6}}
\put(107,90){$a$}
\put(60,75){\circle*{6}}
\put(60,60){$b$}
\put(20,55){\circle*{6}}
\put(15,40){$z'$}
\put(-69,75){\circle*{6}}
\put(-79,65){$c$}
\put(-13,75){\circle*{6}}
\put(-20,80){$z$}
\put(-45,69){\circle*{6}}
\put(-55,55){$d$}
\put(61,220) {\line(1,-3){41}}
\put(60,220) {\line(0,-1){144}}
\put(60,220) {\line(-1,-2){74}}
\put(-69,75) {\line(1,0){130}}
\put(72,188){\line(-2,-5){54}}
\put(72,188){\line(-1,-1){120}}
\put(17,54){\line(2,1){85}}
\put(-69,74){\line(5,-1){87}}
\put(-47,69){\line(5,1){147}}
\end{picture}
\vspace{-10mm}
\caption{}
\end{center}
\end{figure}
 \begin{proof}
Let $\Pa= P_1\uplus \dotsm\uplus P_r$. 
 We claim that $P_r$ is a simplicial vertex of $G_\Pa$. 
Suppose that $\{P_r, P_i\}$ and $\{P_r, P_j\},\, i<j,$ are two
different edges of
$G_\Pa$ and that there are two nontrivial lines $L_1:=\{x,y,z\}$ and 
$L_2=\{x',y',z'\}$ where $x,y, x',y'\in P_r$ and $z\in
P_i,$ 
$z'\in P_j$.
We will consider two possible cases:
\begin{namelist}{xxx}
\item[~$\circ$] 
Suppose
first that two of the elements $x,y, x',y'$ are equal; 
w.l.o.g., we can suppose $x=x'.$ As $\M$ is binary the elements
$x, y, y'$  can't be  colinear, so $\cl(\{x, y, y'\})$ is a plane. 
From modularity of $F_{r-1}$, we know that $\cl(\{x, y, y'\})
\cap F_{r-1}$ is a line. So the line
$\cl(\{y, y'\})$ meets  the modular hyperplane
$F_{r-1}$ in a point $a$. Now the line $\{z,z',a\}$ is a nontrivial line 
which meets $P_i$ and 
$P_j$. Then by definition $\{P_i, P_j\}$ is an edge of $G_\Pa$.
\item[~$\circ$] 
Suppose now that  the elements $x,y, x',y'$ are different.
 Then as $\M$ is binary we have 
$\rk(\{x,y, x', y'\})=4$.
From modularity of $F_{r-1}$,  we know that $\rk(\cl(\{x,y,
x', y'\})\cap F_{r-1})=3$. Then the  six lines $\cl(\{\alpha, \beta\}),$ 
for $\alpha$ and $\beta$ in $\{x,y, x', y'\}$ meet $F_{r-1}$ in
six coplanar points; let these points be labelled
as in  Figure 2. 
Let $P_\ell$ be the set that contains $a$. We
will consider three subcases.
\begin{namelist}{xxx}
\item[~\textbullet] 
Suppose first that $i<j<\ell$. From the property given immediately after
Definition \ref{Pa}, we have that $c$ is also in $P_\ell$.
Consider the modular flat $F_{\ell-1}=\bigcup_{h=1}^{\ell-1} P_h$.
We know that the plane $\cl (\{a, c, z, z'\})$ meets $F_{\ell-1}$ in a line,
so $\cl (\{z,z'\})$ is a nontrivial line meeting $P_i$ and 
$P_j$ and so $\{P_i, P_j\}$ is an edge of $G_\Pa$.
\item[~\textbullet]
Suppose now that $\ell<i<j$. Then the nontrivial line 
$\{a, d, z\}$ meets $P_i$ and $P_\ell$ and
 we have $d\in P_i$. So 
the nontrivial line $\{c, d, z'\}$ meets $P_i$ and $P_j$ and
$\{P_i, P_j\}$ is an edge of $G_\Pa$.
\item[~\textbullet]
Suppose finally that $i\leq \ell\leq j$. 
The nontrivial line $\{a, d, z\}$ meets $P_i$ and $P_\ell$ so 
$d\in P_\ell$.
The nontrivial line $\{c, d, z'\}$ meets $P_\ell$ and $P_j$ and 
necessarily  we have  $c\in P_j$. We conclude that the nontrivial 
line $\{b, c,z\}$ meets  $P_i$
and $P_j$ and $\{P_i, P_j\}$ is an edge of $G_\Pa$.
\end{namelist}
\end{namelist}
By induction we conclude that
$G_\Pa$ is chordal.
\end{proof}
We say that two $M$-chains 
$$\F:=\,\emptyset\subsetneq \dotsm\subsetneq F_{r-1}\subsetneq 
F_r=[n]$$ and
$$\F':=\, \emptyset\subsetneq \dotsm\subsetneq
F'_{r-1}\subsetneq  F'_r=[n]$$
are related by an {\em elementary deformation} if
they differ by at most one flat.
We say that two $M$-chains are {\em equivalent} if they can 
be obtained from each other by   elementary deformations.
\begin{proposition}\label{equiv}
Every two $M$-chains of the same matroid $\M$  are equivalent.
\end{proposition}
\begin{proof} We prove it by induction on the rank. The result is clear 
for $r=2.$ Suppose it true for all matroids of rank at
most $r-1.$ Consider two different $M$-chains 
 $$\F:= \hspace{3mm}\emptyset \subsetneq \dotsm\subsetneq
F_{r-1}\subsetneq  F_r=[n]$$
$$\F':= \hspace{3mm}\emptyset\subsetneq \dotsm\subsetneq
F'_{r-1}\subsetneq  F'_r=[n].$$
Let $F_{\ell}$ be the flat of highest rank  of  the $M$-chain $\F$
contained in $F'_{r-1}$. 
We know that $F_j\cap F_{r-1}^{'}$ ,   $j=0, 1,\dotsc, r$,
is a  modular flat of the matroid $\M$  and that
$$\rk(F_j\cap
F_{r-1}^{'})=j-1,\,\,\,
\text{for}\,\,\, j=\ell+2,\dotsc, r-1.$$ 
Let $\F _0:=\F$ and for $i$ from $1$ to $r-1-\ell$, 
let $\F _i$ be the $M$-chain
$$\emptyset \subsetneq \dotsm\subsetneq
F_{l}\subsetneq F_{\ell+2}\cap F'_{r-1}\subsetneq\dotsm 
F_{\ell+i+1}\cap F'_{r-1}
\subsetneq F_{\ell+i+1}\dotsm\subsetneq [n].$$
We have clearly by definition that for $i$ from $0$ to $r-2-\ell$, 
the  $M$-chains
$\F _i$ and $\F _{i+1}$ are equivalent. This sequence of equivalences shows
that $\F$ is equivalent to  
$\F _{r-1-\ell}.$
Finally by the induction hypothesis we have that $\F'$ is equivalent to 
$\F _{r-1-\ell}$ which concludes the proof.
 
\end{proof}
\begin{remark}
{\em  Proposition ~\ref{equiv} can be used to obtain all the S-labels  of a 
given chordal graph $G$ from a fixed one. If
$G$ is doubly-connected  the number of $M$-chains of $\M(G)$ is equal  
to  the half of the number of such labelings, see
\cite{Stan2}, Proposition 2.8.}
\end{remark}
It is natural to ask if, given a chordal graph $G,$  there is a 
supersolvable matroid
$\M$ together with an $M$-partition $\Pa$ such that $G=G_{\Pa}.$ 
Can the matroid $\M$  be supposed graphic?
The next proposition gives a positive answer to these questions:
\begin{proposition}
Let $G=(V,E)$ be a chordal graph with an\, S-labeling\, $v_1,\dotsc,v_m$
of its vertices, and $\widetilde G$ the extension of \,$G$\, by a vertex 
$v_0$ adjacent to all the vertices, i.e.: $$V_{\widetilde G}=V_G
\cup v_0\hspace{0.5cm}\text{and}\hspace{0.5cm}E_{\widetilde G}=E_G
\cup \{\{v_i,v_0\},\,i=1,\dots, m\}.$$
Then $G_{\widetilde \Pa},$ the derived S-graph of\, $\widetilde G$ for 
the S-labeling $v_0, v_1,\dotsc,v_m$ is isomorphic to \,$G.$
\end{proposition}
\begin{proof}
As $v_0$ is adjacent to every vertex $v_i,\, i=1,\dots, m,$
 it is clear that $v_0, v_1,\dotsc,v_m$ is an S-labeling of $\widetilde G.$
Let\, $\Pa$ and \,$\widetilde \Pa$ denote   the corresponding 
$M$-partitions of the graphic matroids $\M(G)$
and $\M(\widetilde G).$ We have $\Pa=P_1\uplus\cdots \uplus P_{m-1}$ and
$\,\widetilde \Pa= \widetilde P_1(=\{v_0, v_1\}),\uplus \widetilde 
P_2\uplus\cdots \uplus\widetilde P_{m}$ with
$\widetilde P_i=P_{i-1}\cup \{v_o,v_i\},$ for  $ i=2,\dots, m.$ 
Now we can see that if
$\{v_i,v_j\},$
$0\leq i<j\leq m-1,$ is an edge of $G$ then $\{\widetilde P_i, 
\widetilde {P_j}\}$
is an edge of $G_{\widetilde \Pa}.$ From Remark~\ref{sub} 
we get that reciprocally $G_{\widetilde \Pa}$
is a subgraph of $G.$
\end{proof}
\centerline{\textbf {Acknowledgements}}
\vspace{1ex}
\noindent The authors are grateful to the anonymous referees  
for their detailed remarks and suggestions on a previous version of this
paper.

\end{document}